\documentclass{amsart}
\usepackage{amssymb,latexsym}
\usepackage{graphicx}

\newtheorem{proposition}{Proposition}

\def\2{\u a}
\def\3{\c s}
\def\4{\^ a}
\def\5{\c t}
\begin{document}
\thispagestyle{empty}
\phantom{0}
\title{Some properties of regular integers modulo $n$}
\author{Br\2du\5 Apostol}
\begin{abstract}

Denoting by $V(n)$ the number of regular integers $(\text{mod
}n)$, we study properties of the sequences $\displaystyle
\biggl(\frac{V(n+1)}{V(n)}\biggr)_{n\geq1}$ and $\displaystyle
(V(n+1)-V(n))_{n\geq1}$. We also prove that the sets
$\displaystyle\biggl \lbrace\frac{\psi(n)}{V(n)} : n\geq1\biggr
\rbrace$ and $\displaystyle\biggl \lbrace\frac{V(n)}{\phi(n)} :
n\geq1\biggr \rbrace$ are everywhere dense in $(1,+\infty)$.
\end{abstract}

\maketitle
\emph{2010 Mathematics Subject Classification: 11A25, 11N37}

\emph{Key words and phrases:} regular integers $(\text{mod }n)$,
inequalities, density problems.
\section{Introduction}\label{S: intro}

An integer $a$ is called regular\text{ }$(\text{mod }n)$ if there is an integer  $x$ such that $a^{2}x\equiv a\text{ } (\text{mod }n)$

An element $a$ of a ring $R$ is said to be regular (in the sense of J. von Neumann) if there is $x\in R$ such that $a=axa$.
If $R$ is the ring $\mathbb{Z}_n$ of residue classes $(\text{mod }n)$, this condition is equivalent to the congruence above.
Let $\displaystyle \text{Reg}_n=\lbrace{a: 1\leq a \leq n \text{ and  a is regular}\text{ }(\text{mod }n)}\rbrace$
and let $V(n)=\#\text{Reg}_n$.

We notice that $V(p^{\alpha})=p^{\alpha}-p^{{\alpha}-1}+1=\phi(p^{\alpha})+1$, where $\phi$ is the Euler function. The function $V(n)$ is multiplicative and $V(n)=((\phi(p_1^{\alpha_1})+1)...(\phi(p_k^{\alpha_k})+1)$, where $n=p_1^{\alpha_1}...p_k^{\alpha_k}>1$, see ~\cite{oA07}, ~\cite{lT08}.

We denote as usual by $\sigma(n)$ the sum of natural divisors of
$n$ and by $\psi(n)$ the Dedekind function. The formulae
$\sigma(n)={\displaystyle\prod_{i=1}^k{\frac{p_i^{\alpha_i+1}-1}{p_i-1}}}$
and $\psi(n)=n{\displaystyle\prod_{p\vert
n}\biggl(1+\frac{1}{p}\biggr)}$ are well known.
\section{The sequences $\displaystyle\biggl (\frac{V(n+1)}{V(n)}\biggr)_{n\geq 1}$ and $(V(n+1)-V(n))_{n\geq1}$}
\label{S: V(n+1)/V(n) and V(n+1)-V(n)}
In this section  we prove some properties of the sequences $\displaystyle \biggl(\frac{V(n+1)}{V(n)}\biggr)_{n\geq 1}$ and
$\displaystyle \bigl(V(n+1)-V(n)\bigr)_{n\geq1}$.\label{S: Sequences}
We notice that
$V(1)<V(2)$, $V(3)=V(4)$, $V(7)>V(8)$, $V(10)<V(11)$, $V(11)>V(12)$... .
More generally, we prove
\begin{proposition}\label{P: 1}
$\displaystyle\frac{V(n+1)}{V(n)}> 1$ for infinitely many $n$ and, also,
$\displaystyle\frac{V(n+1)}{V(n)}< 1$ for infinitely many $n$.
\end{proposition}
\begin{flushleft}
\textbf{\emph{Proof.}}
Let $p_1,p_2,\ldots,p_r$ be prime numbers. Since
$(1,p_1...p_r)=1$, according to the theorem of Dirichlet, the
arithmetic progression
$1$, $1+p_1...p_r$, $1+2p_1...p_r$, ... contains infinitely many primes.\\
Consider a prime $p$ of the form $p=1+ap_1...p_r$. Taking into account that $V(n)\leq n$, for any $n\geq 1$ and
$V(p)=p$, for any prime $p$,
we have
$$\frac{V(p-1)}{V(p)}=\frac{V(ap_1...p_r)}{1+ap_1...p_r}\leq \frac{ap_1...p_r}{ap_1...p_r+1}<1.$$
Therefore, $\displaystyle \frac{V(n+1)}{V(n)}> 1$
for infinitely many positive integers $n$.
On the other hand, since
$\displaystyle (-1,p_1^2p_2...p_r)=1$, $(p_1,...p_r$ being prime numbers), there are infinitely many primes in the arithmetic progression $-1, -1+p_1^2p_2...p_r, -1+2p_1^2p_2...p_r, ...$ .
Let $q$ be a prime number of the form $q=-1+bp_1^2p_2...p_r$.
It is easy to see that $V(mn)\leq mV(n)$ for every positive integers $m, n\geq 1$.
Therefore
$$\displaystyle \frac{V(q)}{V(q+1)}=\frac{bp_1^2p_2...p_r-1}{V(bp_1^2p_2...p_r)}\geq
\frac{bp_1^2p_2...p_r-1}{b(p_1^2-p_1+1)p_2...p_r}>1.$$
It follows that $\displaystyle \frac{V(n+1)}{V(n)}< 1$
for an infinitely many numbers $n$.
\qquad $\square$\\[10pt]
\end{flushleft}
If $\mathcal{A}=\biggl \lbrace n\in \mathbb{N}: \displaystyle
\frac{V(n+1)}{V(n)}> 1\biggr \rbrace$
and
$\mathcal{B}=\biggl \lbrace n\in \mathbb{N}: \displaystyle
\frac{V(n+1)}{V(n)}< 1\biggr \rbrace$,
another interesting property of  $\displaystyle \frac{V(n+1)}{V(n)}$ is given by\\
\begin{proposition}\label{P: 2}
$\displaystyle \liminf_{\substack{n \to \infty\\n \in
\mathcal{A}}}\frac{V(n+1)}{V(n)}=1$
and
$\displaystyle \limsup_{\substack{n \to \infty\\n \in
\mathcal{B}}}\frac{V(n+1)}{V(n)}=1$.
\end{proposition}
\begin{flushleft}
\textbf{\emph{Proof.}}  We will use Linnik's theorem,
wich states that if $(k,l)=1$, then there is a prime $p$ such that
$p\equiv l \text{ }(\text{mod }k)$ and $p\ll k^c$,
where $c$ is a constant (one can take $c\leq 11$).
Let
$\displaystyle A=\prod_{\substack{p\leq x\\p\text{ prime}}}p=2p_1p_2...p_t$
with $p_1\geq 3$ and $p_1<p_2<...<p_t$ being prime numbers.
Since
$(A^2,A+1)=1$, by Linnik's theorem,
there is a prime number $q$ such that
$q\equiv A+1\text{ }(\text{mod }A^2)$ and $q\ll \bigl(A^2\bigr)^c$= $A^{2c}$,
where $c$ satisfies $c\leq 11$.
Let $q$ be the least prime number satisfying the conditions above.
We have $q-A-1=kA^2$, for some $k$.
Then
$q-1=AB$, where $B=1+kA$. Since $(A,B)=1$,
it follows that $B$ is free of prime factors $\leq x$,
and
$AB=2p_1...p_t\cdot q_1^{\alpha_1}...q_k^{\alpha_k}$, $p_1<...<p_t<q_1<...<q_k$, $p_1\geq 3$.\\
It is known that $\phi(n)<V(n)\leq n$ for every $n> 1$, and
$V(n)=n$ if and only if $n$ is squarefree, see ~\cite{oA07}.
We thus have
$$\displaystyle \frac{V(q)}{V(q-1)}=\frac{V(AB+1)}{V(AB)}=\frac{AB+1}{V(A)V(B)}\geq \frac{AB+1}{AB}>1,$$ so  $q\in \mathcal{A}$.
Since $V$ is multiplicative, we obtain
\begin{equation}\label{Eq: 1}
\tag{$1$}\frac{V({q})}{V(q-1)}=\frac{AB+1}{AB}\cdot\frac{A}{V(A)}\cdot\frac{B}{V(B)}.
\end{equation}
Here
$\displaystyle \frac{AB+1}{AB}\to 1$ as $x \to \infty$,
so it is sufficient to study
$\displaystyle \frac{A}{V(A)}$ and $\displaystyle \frac{B}{V(B)}$.
We have
\begin{equation}\label{Eq: 2}
\tag{$2$}  \frac{A}{V(A)}=\frac{2p_1p_2...p_t}{V(2p_1p_2...p_t)}=1.
\end{equation}
It is well-known that
$A=\displaystyle \prod_{p\leq x}p=e^{O(x)}$.
Since $q\ll A^{2c}$, $A=e^{O(x)}$ and $B\ll A^{10}$, we obtain $B\ll \bigl( e^{O(x)}\bigr)^{10}= e^{O(x)}$, so
\begin{equation}\label{Eq: 3}
\tag{$3$}  \log B\ll x.
\end{equation}
If $B=\displaystyle \prod_{i=1}^{k}q_i^{b_i}$ is the prime factorization of $B$, then\\
$\log B=\displaystyle \sum_{i=1}^{k}b_i\log q_i>(\log x)\sum_{i=1}^{k}b_i$,
as $q_i>x$ for all $i\in \lbrace1,2,...,k \rbrace$. Here
$\displaystyle \sum_{i=1}^{k}b_i\geq k$, thus $\log B>k\log x$.
It follows by \eqref{Eq: 3} that
\begin{equation}\label{Eq: 4}
\tag{$4$}   k<\frac{\log B}{\log x}\ll \frac{x}{\log x}.
\end{equation}
For $\displaystyle \frac{B}{V(B)}$, we have
$\displaystyle \frac{B}{V(B)}=\prod_{i=1}^k\frac{q_i^{b_i}}{q_i^{b_i}-q_i^{b_i-1}+1}<
\prod_{i=1}^k\frac{q_i^{b_i}}{q_i^{b_i}-q_i^{b_i-1}}=\prod_{i=1}^k\frac{1}{1-\frac{1}{q_i}}$.
Using $\displaystyle 1-\frac{1}{q_i}>1-\frac{1}{x}$ and \eqref{Eq: 4}
we obtain
$$\displaystyle \prod_{i=1}^{k}\biggl(1-\frac{1}{q_i} \biggr)>\biggl(1-\frac{1}{x}\biggr)^k\geq
\biggl(1-\frac{1}{x}\biggr)^{O(\frac{x}{\log x})}=1+O\biggl (\frac{1}{\log x} \biggr),$$
so,
\begin{equation}\label{Eq: 5}
\tag{$5$}   \frac{B}{V(B)}<\frac{1}{1+O\bigl (\frac{1}{\log x} \bigr)}.
\end{equation}
By \eqref{Eq: 1}, \eqref{Eq: 2}, \eqref{Eq: 5} and
$\displaystyle \frac{AB+1}{AB}\to 1$ as $x\to \infty$, it follows that\\
\begin{equation}\label{Eq: 6}
\tag{$6$}   \frac{V(q)}{V(q-1)}<\frac{1}{1+O\bigl (\frac{1}{\log x} \bigr)},
\end{equation}
By \eqref{Eq: 6} from above and $q\in \mathcal{A}$,
it follows that
$\displaystyle \liminf_{\substack{n \to \infty\\n \in \mathcal{A}}}\frac{V(n+1)}{V(n)}=1$,
so the first part of the proposition is proved.\\
Now let
$\displaystyle A=p_{\pi(x)}\prod_{\substack{p\leq x\\p\text{
prime}}}p=2p_1...p_{t-1}p_{t}^2$ , with $p_1\geq 3$ and
$p_1<p_2<...<p_t$ being prime numbers ($t=\pi(x)$, the number of prime numbers $\leq x$, where $x>0$).
Since $(A^2,A-1)=1$, by Linnik's theorem,
there is a prime number $q$ such that $q\equiv A-1\text{ }(\text{mod }A^2)$.
Let $q$ be the least prime number satisfyng conditions from above,
so $q=AB-1$, where $B=1+kA$ is free of prime factors $\leq x$.
If $n$ is not a squarefree number, it is easy to see that $V(n)\leq n-2$, for sufficiently large $n$.
Take $n=q$.
Since $AB$ is not a squarefree number, we have
$$\frac{V(q+1)}{V(q)}=\frac{V(AB)}{AB-1}\leq \frac{AB-2}{AB-1}<1,$$
so $q\in \mathcal{B}$.
We have
$$\displaystyle \frac{V(q+1)}{V(q)}=\frac{V(AB)}{AB-1}=
\frac{V(A)}{A}\cdot \frac{V(B)}{B}\cdot \frac{AB}{AB-1}.$$
Then,
$\displaystyle \frac{V(A)}{A}=\frac{p_t^2-p_t+1}{p_t^2}$ and
$$\displaystyle A=p_{\pi(x)}\prod_{\substack{p\leq x\\p\text{ prime}}}p=O(x)e^{O(x)}=e^{\log O(x) +O(x)}=
e^{O(x)}.$$
Using a similar argument as in the first part, we obtain that
$$\displaystyle \frac{V(B)}{B}>1+O\biggl (\frac{1}{\log x} \biggr).$$
Since $\displaystyle \frac{AB}{AB-1} \to 1$ as $x  \to \infty$,
$\displaystyle \frac{p_t^2-p_t+1}{p_t^2}\to 1$ for $x  \to \infty$
and $\displaystyle \frac{V(B)}{B}>1+O\biggl (\frac{1}{\log x}
\biggr)$, we deduce that $\displaystyle \limsup_{\substack{n \to
\infty\\n \in \mathcal{B}}}\frac{V(n+1)}{V(n)}=1$. \qquad
$\square$
\end{flushleft}
In connection with the difference $V(n+1)-V(n)$, we can prove
\begin{proposition}\label{P: 3}
$$\limsup_{n \to \infty}(V(n+1)-V(n))=+\infty,\text{ } and\text{ }
\liminf_{n \to \infty}(V(n+1)-V(n))=-\infty.$$
\end{proposition}
\begin{flushleft}
\textbf{\emph{Proof.}}
Let $n=2^km$, ($k\geq 1$, $m$ odd) be an even number. Since $V$ is multiplicative,
$V(n)=V(m)V(2^k)\leq m(2^k-2^{k-1}+1)=\displaystyle \frac{n}{2}+\frac{n}{2^k}$.
So,
$\displaystyle V(n)\leq \frac{3}{4}n$ for every $n$ which is a multiple of 4.\\
Let $p$ be a prime number of the form $p=4t+1$. Then $V(p)=p$,
so
$$\displaystyle V(p)-V(p-1)\geq p-\frac{3}{4}(p-1)=\frac{p+3}{4}.$$
Since (according to the theorem of Dirichlet) we may take $p$ as large as we please,
the first assertion is proved.\\
Now take $p$ a prime number of the form  $p=4t+3$. Then
$$\displaystyle V(p)-V(p+1)\geq p-\frac{3}{4}(p+1)=\frac{p-3}{4}.$$
Since $p$ in the above relation may be taken arbitrarily large,
the second assertion is proved.  \qquad $\square$
\end{flushleft}
\section{Density problems}\label{S: Density}
We study the density of the sets
$\displaystyle\biggl \lbrace\frac{\psi(n)}{V(n)} : n\geq1\biggr \rbrace$
and
$\displaystyle\biggl \lbrace\frac{V(n)}{\phi(n)} : n\geq1\biggr \rbrace$.
In order to prove the next proposition we apply the following
result of B.S.K.R. Somayajulu, ~\cite{bS77}: Let $(a_n)_{n\geq 1}$
be a strictly decreasing sequence of positive numbers with
$\displaystyle \lim_{n \to\infty}{a_n}=0$
(we use the notation $0<a_n\searrow 0$), and suppose that the series $\displaystyle \sum{a_n}$ is divergent.
Then, for each $s>0$, there is an infinite subseries of  $\displaystyle \sum{a_n}$, which converges to $s$.
\begin{proposition}\label{P: 4}
\begin{equation}
\tag{$i$}  The\text{ } set\text{ } \biggl \lbrace\frac{\psi(n)}{V(n)} : n\geq1\biggr \rbrace\text{ } is\text{ } everywhere
\text{ } dense\text{ } in\text{ }
(1,+\infty),
\end{equation}
\begin{equation}
\tag{$ii$}  The\text{ } set\text{ } \biggl \lbrace\frac{V(n)}{\sigma(n)} : n\geq1\biggr \rbrace\text{ } is\text{ } everywhere
\text{ } dense\text{ } in\text{ }
(1,+\infty).
\end{equation}
\end{proposition}
\begin{flushleft}
\textbf{\emph{Proof.}}
\end{flushleft}
\begin{flushleft}
\textbf{$(i)$}\quad  Let $p_n$ the $n-$th prime and  $\displaystyle u_n=\frac{1}{p_n}$.
Then, $0<u_n\searrow 0$ for $n \to \infty$ and
$$\displaystyle \prod_{n=1}^{\infty}(1+u_n)=\prod_{n=1}^{\infty}(1+\frac{1}{p_n})=\infty.$$
Now let $a_n=\log(1+u_n)$. We have $0<a_n\searrow 0$, for $n \to \infty$.
Then
$$\displaystyle \sum_{n=1}^{\infty}a_n= \sum_{n=1}^{\infty}\log(1+u_n)=\log\prod_{n=1}^{\infty}(1+u_n)=\infty,$$
so the series
$\displaystyle \sum_{n=1}^{\infty}a_n$ diverges.
Taking into account the result of Somayajulu, for each $\delta>1$ there exists $(a_{{n}_{r}})_{r\geq 1}$, such that that the subseries
$\displaystyle \sum a_{{n}_{r}}$ converges and $\displaystyle \sum_{r=1}^{\infty}a_{{n}_{r}}=\log{\delta}$.
We have
$$\displaystyle \sum_{r=1}^{\infty} a_{{n}_{r}}= \sum_{r=1}^{\infty}\log(1+u_{{n}_{r}})=
\log\prod_{r=1}^{\infty}(1+u_{{n}_{r}}),$$
so
$\displaystyle \prod_{r=1}^{\infty}(1+u_{{n}_{r}})=\delta$.
Take $m_r=p_{{n}_{1}}...p_{{n}_{r}}$.
Then
$$\displaystyle \frac{\psi(m_r)}{V(m_r)}=\frac{(p_{{n}_{1}}+1)...(p_{{n}_{r}}+1)}{p_{{n}_{1}}...p_{{n}_{r}}}=
\prod_{k=1}^{r}\biggl(1+\frac{1}{p_{{n}_{k}}} \biggr)=\prod_{k=1}^{r}(1+u_{{n}_{k}}),$$
and
$\displaystyle \lim_{r \to \infty}\frac{\psi(m_r)}{V(m_r)}=\lim_{r \to \infty}\prod_{k=1}^{r}(1+u_{{n}_{k}})=\delta$.
So,
$\displaystyle\biggl \lbrace\frac{\psi(n)}{V(n)} : n\geq1\biggr \rbrace$ is everywhere dense in
$(1,+\infty)$.
\end{flushleft}
\begin{flushleft}
\textbf{$(ii)$}\quad  Choose $u_n$ with $1+u_n=\displaystyle \frac{1}{1-\frac{1}{p_n}}$.
Then
$u_n=\displaystyle \frac{1}{1-\frac{1}{p_n}}-1=\frac{1}{p_n-1}>\frac{1}{p_{n+1}-1}=u_{n+1}$ and
$\displaystyle \lim_{n \to \infty}u_n=0$.
So,\\
$0<u_n\searrow 0$ for $n \to \infty$.\\
Since $1+u_i=\displaystyle \frac{1}{1-\frac{1}{p_i}}>1+\frac{1}{p_i}$,
it follows that
$$\displaystyle \prod_{i=1}^n(1+u_i)=\prod_{i=1}^n\frac{1}{1-\frac{1}{p_i}}>
\prod_{i=1}^n\biggl(1+\frac{1}{p_i}\biggr),$$
so
$\displaystyle \prod_{n=1}^{\infty}(1+u_n)=\infty$. Now let
$a_n=\log(1+u_n)$. Then $0<a_n\searrow 0$, for $n \to \infty$.\\
Since $\displaystyle \sum a_n$ diverges, for each $\delta>1$ there exists $(a_{n_{r}})_{r\geq 1}$, such that the subseries
$\displaystyle \sum a_{n_{r}}$ converges and $\displaystyle \sum_{r=1}^{\infty}a_{n_{r}}=\log\delta$.
Moreover,
$\displaystyle \prod_{r=1}^{\infty}(1+u_{n_{r}})=\delta$.\\
Take $m_r=p_{{n}_{1}}...p_{{n}_{r}}$.
Then
$$\displaystyle \frac{V(m_r)}{\phi(m_r)}=\frac{p_{{n}_{1}}...p_{{n}_{r}}}{(p_{{n}_{1}}-1)...(p_{{n}_{r}}-1)}=
\displaystyle \prod_{k=1}^{r}\frac{1}{1-\frac{1}{p_{n_{k}}}}= \prod_{k=1}^{r}(1+u_{n_{k}})\to \delta \text{ } \text{as}
\text{ } k\to \infty.$$
This completes the proof of $(ii)$. \qquad $\square$
\end{flushleft}

ACKNOWLEDGEMENT. The author would like to thank  Professor L.T$\acute{o}$th and Professor G. Mincu for the 
valuable suggestions which improved the presentation.

\textbf{Br\2du\5 Apostol}, "Spiru Haret" Pedagogic High School, 5 Timotei Cipariu St., RO-620004  Foc\3ani, Rom\4nia,
 E-mail: apo\_brad@yahoo.com\\[8pt]
\end{document}